\documentclass[11pt]{article}
\usepackage{amssymb,amsmath,graphicx}   

\newcommand{\be}{\begin{eqnarray}}
\newcommand{\nd}{\end{eqnarray}}

\begin{document}
\begin{center}{\bf On resolution to Wu's Conjecture}\\
by Theodore Yaotsu Wu\footnote{California Institute of Technology, Pasadena, CA ~91125 U.S.A. ~Email:tywu@caltech.edu} \end{center}
\noindent{\bf Abstract.}~In this series of studies on Cauchy's function $f(z)$ ($z=x+iy$) and its integral $J[f(z)]\equiv (2\pi i)^{-1}\oint_C f(t)dt/(t-z)$ taken along a Jordan contour $C$, the aim is to investigate their comprehensive properties over the entire $z$-plane consisted of the simply-connected closed domain ${\cal D}^+$ bounded by $C$ and the open domain ${\cal D}^-$ outside $C$. This article attempts to solve an inverse problem that Cauchy function $f(z)$, regular in ${\cal D^+}$ and on $C$, has a singularity distribution in  ${\cal D}^-$ which can be determined in analytical form in terms of the values $f(t)$ numerically prescribed on $C$, which is Wu's conjecture[1].  It is resolved here for $f(z)$ having (i) a single, (ii) double, or (iii) multiple singularities of the types (I) $\Sigma_j^N M_j(z_j-z)^{k_j}$, (II) $M_\ell\log(z_2-z)$, by having their power series expanded in $z$ and matched on a unit circle $(t=e^{i\theta}, -\pi\leq\theta<\pi$ for contour $C$) with the numerically prescribed Fourier series $f(z)=\Sigma_0^\infty c_ne^{in\theta}$ for solution.  The mathematical methods used include (a) complex algebra for cases (i)-(ii), (b) for case (iii) a general asymptotic method developed here for resolution to the Conjecture by induction, and (c) the generalized Hilbert transforms to expound essential singularities.  This Conjecture has an advanced version for $f(z)$ to be given only one of its two conjugate functions on $C$ to suffice, and another for the complement function $F(z)$ defined as being regular in domain ${\cal D}^-$ and having singularities in ${\cal D^+}$.  These new methods are applicable to all relevant problems in mathematics, engineering and mathematical physics requiring breakthrough by having the exterior singularities resolved.

\vskip 1.6mm
\noindent{\bf 1. Introduction.}
\vskip 0.6mm
The foundation of the present study starts with Cauchy's theorem and integral formula,
\vskip -2mm
\begin{subequations}\label{1}
\be   J[f(z)]\equiv\frac{1}{2\pi i}\oint_C \frac{f(t)}{t-z}dt &=& f(z) \qquad~~~ (z\in {\cal D^+} -\mbox{open domain inside} ~C), \label{1a}\\
    &=& 0   \qquad\qquad~ (z\in {\cal D^-} -\mbox{open domain outside} ~C), \label{1b} \nd
\end{subequations}
where Cauchy's function $f(z)$ is assumed to be analytic, regular $\forall z\in {\cal D^+}$ and continuous for $z=t$ on contour $C$, taken in positive (counter-clockwise) sense as understood.  Here, (\ref{1b}) follows from {\it Cauchy's integral theorem} that $\oint_C g(t)dt=0$ if $g(t)$ is regular within and on contour $C$ (as is $g(t)=f(t)/(t-z) ~\forall z\in {\cal D^-}$), while (\ref{1a}) is known as {\it Cauchy's integral formula}, also called {\it Cauchy's functional relation}, known to hold for $z$ in open domain ${\cal D^+}$ (but not including its boundary contour $C$).
\vskip 0.8mm
The task of determining the value of Cauchy's integral $J[f(z)]$ for $z$ situated right on $C$ has been accomplished by Wu[1] with adopting a generalized condition that
\be  f(z)~\mbox{be}~C^{n}~\forall z\in {\cal D^+}~\mbox{and in a neighborhood ${\cal N}_C$ striding across contour}~C ~~(n~\mbox{being arbitrary}),  \label{2}\nd
where the corresponding function $f(z)$ is called the {\it generalized Cauchy's function}.  This new condition thus renders the integral $J[f(z)]\equiv (2\pi i)^{-1}\oint_C f(t)dt/(t-z)$ intact in value as a point $z\in {\cal D^+}$ (or $z\in {\cal D^-}$) tends to a generic point $z_o$ on $C$ without crossing the contour if the contour remains fixed except for an infinitesimal stretch about $z_o$ where it is indented into a semi-circle $C^+_{\epsilon}$ onto the ${\cal D^-}$-side (or a semi-circle $C^-_{\epsilon}$ onto the ${\cal D^+}$-side), of radius $\epsilon$ from $z_o$ on $C$ (see Wu[1], Fig. 1). In the limit as $\epsilon\rightarrow 0$, the integral $J[f(z)]$ having $z\{\in {\cal D^+}\}\rightarrow z_o$ (from inside $C$) is let to assume its inner limit of undetermined value $f^+(z_o)$ (resulting from the integral over the semi-circle $C^+_{\epsilon}$ plus the integral over $C-C^+_{\epsilon}$ which assumes its principal value as $\epsilon \rightarrow 0$).  Similarly, the integral $J[f(z)]$ having $z\{\in {\cal D^-}\}\rightarrow z_o$ (from outside $C$) has its outer limit of $J[f(z_o)]=f^-(z)=0$ as is invoked by Cauchy's Theorem (\ref{1b}), and this limit is the sum of the integral over the semi-circle $C^-_{\epsilon}$ plus the principal value of the integral over $C-C^-_{\epsilon}$.  The final two limit equations thus yield three key relations (cf. Wu[1], Eq. (12)) as
\be (I): f^+(z)=f(z); \quad~~(II):f^-(z)=0;\quad~~(III): f(z)=\frac{1}{\pi i}{\cal P}\oint_{C}\frac{f(t)}{t-z}dt \qquad (z ~\in~ C), \label{3} \nd
in which the integral with symbol ${\cal P}$ denotes its Cauchy principal value while the suffix of $z_0$ is omitted for all $z$ on $C$.  Here, relation (I), $f^+(z)=f(z)$, shows that the limit $f^+(z)$ of integral $J[f(z)]$ resulting as $z~\{\in {\cal D^+}\}\rightarrow z_o$ is equal to the original prescribed $f(z) ~\forall z\in C$, therefore proves the {\it uniform continuity of $f(z)$ in the closed domain} $\overline{{\cal D}^+}=[{\cal D}^+ +C]$, relation (II) shows $J[f(z)]=0 ~\forall~z$ in the closed domain $\overline{{\cal D}^-}=[{\cal D}^- +C]$, whilst relation (III) relates $f(z)$ for each $z\in C$ in terms of all the other values of $f(t)$ over $C$.  Since the original Cauchy's two integral theorems in (\ref{1}) are thus shown also holding valid for $z\in C$, this further establishes the theorem that {\it Cauchy's integral $J[f(z)]$ is uniformly convergent in the closed domain} $\overline{{\cal D}^+}=[{\cal D}^+ +C]$.  These are among the important results discovered in Wu[1].
\vskip 0.8mm
The key relations (I)-(III) in (\ref{3}) have various prospects for further applications and development.  When the contour $C$, originally arbitrary in shape, assumes certain geometric forms in particular, e.g. one that circumventing the upper-half (or lower-half) of the $z$-plane, or another inside (or outside) the unit circle $|z|=1$, the specific Hilbert transform relations between the conjugate functions $u$ and $v$ of $f(z)=u(x,y)+iv(x,y)$ follow directly from the formulas in (\ref{3}) (cf. Wu[1]), {\it which can therefore be regarded as providing the generalized Hilbert transform formulas with a simple proof}.
\vskip 0.8mm
The relations (I)-(III) further help illustrate the fact that a unique relation exists between the values of an analytic function $f(t)~\forall z=t $  varying on contour $C$ and all the singularities of $f(z)~\forall z \in {\cal D^-}$ outside $C$, as exemplified in explicit cases of the so-called {\it direct problems defined by having the singularities prescribed} (cf. Wu[1], \S 7).  The only exception to this assertion of ever existence of singularities of $f(z)~\forall z \in {\cal D^-}$ is when $f(t)\equiv 1$ (or any constant) on $C$, in which case $f(z)\equiv 1 ~\forall z$ in the entire $z$-plane, including $z=\infty$, as can be shown by analytic continuation, in virtue of Liouville's theorem.
\vskip 0.8mm
The present study is motivated by Wu's[1] {\it inverse problem} asserting that every Cauchy function $f(z)$ invariably has a singularity distribution in domain ${\cal D}^-$ outside $C$ which can be determined in closed analytical form in terms of the values $f(t)$ prescribed only numerically on contour $C$, i.e.
\vskip 0.8mm
\noindent{\bf The inverse problem.} ~{\it The inverse problem is (i) to have function $f(t)$ prescribed over contour $C$ only numerically, or equivalently in terms of a series with known numerical coefficients, (ii) $f(z)$ is regular inside contour $C$, and (iii) to use the given numerical data to determine the entire singularity distribution of $f(z)~\forall z\in \cal D^-$ outside $C$ in a closed analytical form, whatever the singularity distribution.}  This seems to cover such unsettled situations as resulting from using, e.g. perturbation series expansion, boundary integral scheme, or other computational algorithms.  For the general solution to the inverse problem, its existence is enunciated in terms of the following Conjecture (cf. Wu[1], \S 10.3):
\vskip 0.8mm
\noindent{\bf The conjecture.} ~{\it A solution to this inverse problem is conjectured to exist, which may not be unique.}\vskip 0.8mm
Resolution to this inverse problem is of vital importance.  This can be well expounded by an expository review of the interesting and colorful history of the theory of water waves, a subject of intimate significance to, and long desired for the general resolution to the Conjecture at hand.
\vskip 0.8mm
Historically, it dates back to the pioneering work of Sir George Gabriel Stokes[2] who developed the perturbation expansion theory (1847) and used it for the very first time to study the nonlinear effects in simple harmonic water waves on deep water to three terms of the power series (expanded in the physical plane) in terms of a small parameter $\alpha=ka, ~a$ being the wave amplitude, $k$ the wave number.  The crucial results attained there were found to produce important qualitative changes in wave behavior and manifest new phenomena, such as the dispersion relation involves the wave amplitude.  Later in 1880, Stokes pursued the expansion for periodic waves on water of constant rest depth $h$, this time formulated in the complex potential plane (for its prefered fixed boundary of the flow region, thereby lending much simpler algebra without involving iterations) to five higher order terms[2].  Along this line, perturbation expansions have been carried out for solitary waves by various authors, including Longuet-Higgins \& Fox[3], Williams[4] and others to further higher orders as summed in review[5], yet still falling short of a conclusion on the series convergence.  Recently, the series expansion for solitary waves of all heights have been brought by Qu[6] to fifteenth order in a series of exact coefficients in rational numbers.  This has been further pursued by Wu et al.[7] to discover that the series is asymptotic in nature within the range of $12-17$th order (varying due to which one or another small parameter adopted for expansion).  But for all the first integrals (e.g. of the excess mass, wave energies, etc), convergence was found by Wang \& Wu[8] to arise within the range of $20-25$th order of the power series expanded only in terms of $\epsilon=k^2h^2$, a final result which is not shared by power series in other parameters, e.g. $a/h$, attested.
\vskip 0.8mm
On the other hand, the issue on the exterior singularities has also been apprehended as a new interest underlying this general problem.  Grant[9] considered the highest periodic water wave (with Stokes's $120^\circ$ corner crest), formulated with the physical variable as a function $z=x+iy=z(f)$ in the lower half of the complex potential $f=\phi+i\psi$ plane continued analytically to the entire $f$-plane.  He showed by analysis that the singularities of $z(f)$ outside the flow field $\overline{{\cal D}^+}$ are all of the regular type of order $1/2$, i.e. $z(f)\sim c(f-f_o)^{1/2}$ as $f\rightarrow f_o~\forall f_o\notin\overline{{\cal D}^+}$.  Grant contended that there exist several singularities of order $1/2 ~\notin\overline{{\cal D}^+}$ which coalesce at the cornered wave crest to become a singularity of order $2/3$ plus a secondary singularity of an irrational order, hence not regular, in a quite perplexing way.  Along a computational approach, Schwartz[10] examined the singularity using the Pad$\acute{e}$ approximation and Domb-Sykes[11] plots, finding that the singularity varies continuously with the wave height from order $1/2$ to $1/3$ (for the complex velocity $df/dz$), and drawing a conjecture that the changes involve coalescence of several singularities of order $1/2$.  However, Tanveer[12] examined the number of singularities, based on the argument principle, finding that there is just one singularity outside the flow field, in disagreement with Schwartz's results.  This question was further investigated by Longuet-Higgins \& Fox[13] with analytic continuation of their $z(f)$ for almost-highest wave expanded in a power series to sixty terms with the Pad$\acute{e}$ approximation to find an exterior branch cut of order $1/2$, in agreement with Grant[9].  However, the conjectured coalescence of singularities into one to fit the cornered crest of the highest wave remains to be expounded.  Finally, any exterior singularity must have its location and magnitude, in addition to its type, also determined to complete the resolution.  This is the main objective of the present study.

Summing up the foregoing two issues, it appears that their eventual resolution may be attained together with first resolving the Conjecture at hand.

We now return to our present endeavor for resolution of Wu's Conjecture.  First, we classify function $f(z)$ to have (i) a single, (ii) double, or (iii) multiple singularities of the algebraic and/or logarithmic types, each of magnitude $M_j$, located at $z_j$, of type $k_j$, to be determined in terms of the numerically prescribed $f(z)$ at the unit circle $z=e^{i\theta}$ as the standard contour $C$.  For cases (i)-(ii), solutions can be obtained using complex algebra for each singularity of $f(z)$ in ${\cal D}^-$ outside $C$, as illustrated in \S 2-\S3.  The original Conjecture is extended to an advanced version for $f(z)$ to have only one of its two conjugate functions prescribed on $C$ to suffice, as shown in \S5, and further for the complement function $F(z)$ defined as being regular in domain ${\cal D}^-$ and having singularities in ${\cal D}^+$.  The case for singularities located on a circle in ${\cal D}^-$ concentric to $C$ is discussed in \S6.  For case (iii), solutions by using complex analysis are shown in \S7, and also by a general method which is developed in \S8 for resolution to the Conjecture by induction.  Finally, the essential singularities are expounded in \S9.
\vskip 1mm
\noindent{\bf 2. Generalized Cauchy's function with a single singularity.}
\vskip 0.8mm
First, we consider the primary case when $f(z)$ (regular within and on contour $C$) possesses a single singularity (of number $N_s=1$) to be determined in terms of its values $f(t)$ prescribed numerically for $z=t~\forall t\in C$, however with no further data available on the base parameters $(M, z_1, k)$ for magnitude $M$, location $z_1$, and {\it primary key index} $k$ (for the order, e.g. $M(z_1-z)^k$) of the sole singularity.  While being primary, this is the basic case to which more general cases can be reduced or referenced.
\vskip 0.8mm
For simplicity, we take the unit circle $|t|=1$ for the standard contour $C$ throughout, onto which any arbitrary regular contour enclosing a simply-connected domain can be mapped by conformal transformation, so that in polar form, $z=r\exp (i\theta), ~r=|z|, ~\theta=\arg z, ~r=1$ on $C$, on which we have
\be  f(t) = f(e^{i\theta})=\tilde{f}(\theta) \qquad~(-\pi\leq\theta\leq \pi), \label{4}\nd
which is prescribed numerically to a desired accuracy so as to satisfy condition (\ref{2}).
\vskip 0.8mm
With this provision, we can then express $\tilde{f}(\theta)$ in a Fourier series in complex form as
\be  \tilde{f}(\theta)=\sum_{n=0}^N ~c_n e^{in\theta}, \qquad~ c_n=\frac{1}{2\pi}\int_{-\pi}^\pi e^{-in\theta}\tilde{f}(\theta)d\theta , \label{5}\nd
with $N$ numerically set for the accuracy desired. This series representation may be regarded as entirely equivalent to $\tilde{f}(\theta)$ given numerically in (\ref{4}), for they are merely a quadrature apart.
\vskip 0.8mm
On the other hand, the sole point singularity can be represented generically as
\begin{subequations}\label{6}
\be  f(z)&=&M(z_1-z)^k \qquad\quad~~ (|z_1|> 1), \label{6a}\\[2mm]
   f_\ell(z)&=& M_\ell \log (z_1-z), \label{6b}\nd\end{subequations}
where $k\in R$, assumed real, so that it is a zero (or a pole) for $k$ being a positive (or a negative) integer, or an algebraic branch connecting $z_1$ and $z=\infty$ for $k$ not an integer, or a logarithmic branch for $k=0$, all lying outside the unit circle for $|z_1|> 1$.  Thus, (\ref{6a})-(\ref{6b}) represent $f(z)$ with a single singularity for resolution, i.e. to have $(M, z_1, k)$ determined in terms of the given coefficients $\{c_n\}$.
\vskip 0.8mm
To proceed, we also expand $f(z)$ of (\ref{6}) into a power series in $z$, yielding
\begin{subequations}\label{7}
\be  f(z)= Mz_1^{k} \sum_{0} \gamma_n\left(-\frac{z}{z_1} \right)^n,~~~\gamma_n(k)=\left(\begin{array}{c}k\\n\end{array}\right)=\frac{k(k-1)\cdots(k-n+1)}{n!}=\frac{\Gamma (k+1)}{\Gamma (k-n+1)n!} , \label{7a}\quad~&&\\[1mm]
  f_\ell(z)=M_\ell \log (z_1-z)=M_\ell \log (z_1)-M_\ell\sum_{n=1}\frac{1}{n}\left(\frac{z}{z_1}\right)^n,\label{7b}\qquad\qquad\qquad\qquad\quad&&\nd\end{subequations}
where $\Gamma(k)$ is the Gamma function, $\Gamma(k+1)=k\Gamma(k), ~\Gamma(n+1)=n!$.
\vskip 0.8mm
Next, matching series (\ref{7a}) termwise with series (\ref{5}) readily yields
\begin{subequations}\label{8}
\be  f(z)= L \sum_{0} \gamma_n\left(-\frac{z}{z_1} \right)^n \equiv s(z),~~~c_o=L\gamma_o,~~~c_n=L\gamma_n (-z_1)^{-n}, \qquad~  (L= Mz_1^{k})\label{8a} \nd
agreeing with (\ref{5}) for $z=e^{i\theta}$ on contour $C$.  This series has the ratio of consecutive coefficients as
\be  R_n(k,z_1)=\frac{c_{n+1}}{c_n}=-\frac{\gamma_{n+1}}{\gamma_n}\frac{1}{z_1}=\frac{n-k}{n+1}\frac{1}{z_1}=\left(1-\frac{1+k}{n+1}\right)\frac{1}{z_1} \quad~~~(n=1,2,\cdots), \label{8b} \nd
which shows (e.g. by the ratio test) that series (\ref{7}) converges for $|z/z_1|<1$ within the circle of convergence of radius $\rho=|z_1|$ on which the point singularity lies at $z=z_1$.  Finally, the undetermined singularity location $z_1$ can be eliminated by taking the ratio $D_n(k)=R_{n+1}/R_n$, giving
\be D_n(k) =\frac{R_{n+1}}{R_n}=\frac{n+1}{n+2} \frac{n+1-k}{n-k}=c_{n+2}~c_n/c_{n+1}^2\quad~~~(n=1,2,\cdots), \label{8c}\nd\end{subequations}
from which $k$ is attained in terms of the given data $c_{n+2}c_n/c_{n+1}^2$ to identify the singularity type.
\vskip 0.8mm
Further, it is clear that the preceding relations (\ref{6})-(\ref{8}) for $f(z)$ also hold valid for the sole logarithmic singularity $f_{\ell}(z)$ of (\ref{6b}) for $k=0$.
\vskip 0.8mm
At last, we note that if (i) $f(z)$ has a single singularity, (ii) the given numerical data $\{c_n\}$'s are exact, then the above relations for $n>3$ all turn out to be identical in effect to the first three orders for $n=1,2,3$ which can provide the final solution.  We therefore have the following theorem proved.
\vskip 0.8mm
\noindent{\bf Theorem 1. Primary resolution to the Conjecture.} ~{\it If (i) function $f(z)$ is regular in domain ${\cal D^+}$ bounded by the unit circle $C~(z=e^{i\theta})$ and on $C$, (ii) $f(z)$ is numerically prescribed on $C$ by series (\ref{5}) with coefficients $\{c_n, ~n=0,1,2,\cdots\}$ which are exact in value, and (iii) $f(z)$ has a single singularity of unknown magnitude $M$, location $z_1~(|z_1|\geq 1)$, and key index $k$ (of order $(z-z_1)^k$), then}
\be  k&=&n+\frac{n+1}{n+1-(n+2)D_n}~(\equiv k(n)), \quad~ D_n=\frac{R_{n+1}}{R_n}, \quad~~R_n=\frac{c_{n+1}}{c_n}, \notag\\[1.2mm]
  z_1&=& \frac{n-k}{(n+1)R_n}~(\equiv z_1(n)), \quad~M=c_0/z_1^k, \quad~ M_\ell =c_0/\log (z_1). \label{9}\nd
\vskip 0.8mm
In principle, the formulas in (\ref{6})-(\ref{9}), if exact, should produce equal values to $(M, z_1, k)$ for every order $n$.  In practice, however, the numerical data for the coefficients $\{c_n\}$'s may contain errors.  Even more so, the presumption for the present case of $N_s=1$ that there is only one singularity needs first to be confirmed.  In this premise, it raises the following important issue to be addressed next.

\vskip 1mm
\noindent{\bf 3. Verification of consistency and error estimate.}
\vskip 0.8mm
If the final results for base parameters $(M, z_1, k)$ obtained from a specific trial run turn out to be uniformly consistent over order $n$ except for very small fluctuations from order to order, which may be due to numerical errors incurred in deriving the coefficients $\{c_n\}$'s for prescribing the given data.  In such cases, the resolution may be said to bear the {\it confirming feature of consistency} in verifying the sole hidden singularity.  If in contrast, the resulting values for $(M,z_1,k)$ appear to vary in finite magnitude for differing orders, such inconsistency then indicates that the correct solution is yet to be found.
\vskip 0.8mm
For error estimate, we assign a small error $\epsilon_n$ to each $c_n$ to give the {\it rounded coefficient} $\breve{c}_n$ and
\be  \breve{R}_n=\breve{c}_{n+1}/\breve{c}_n=(c_{n+1}+\epsilon_{n+1})/(c_n+\epsilon_n),\quad\longrightarrow\quad \breve{k}(n)=k+\kappa(\epsilon_n,\epsilon_{n+1}, \epsilon_{n+2}), \label{10}\nd
so that $\kappa$, the error to $k$, being by (\ref{8c}) an algebraic function of three specific $\epsilon$'s, vanishes with the $\epsilon$'s instead of being nonlinearly magnified.  A similar conclusion can be drawn for ($z_1, M, M_\ell$).
\vskip 0.8mm
\noindent{\bf Example 1.}~To illustrate these points, we let series (\ref{5}) be prescribed for a specific case by 
\be c_0=2,~~~c_1=1.999895,~~~c_2=1.333415,~~~c_3=0.889123,~~~c_4=0.592593,~~~c_5=0.395062, \notag\nd
for the five leading terms taken here for solution and verification. We then find, by (\ref{9}),
\be R_1=c_2/c_1=0.666743,&&R_2=0.666801,~~R_3=0.666492,~~R_4=0.666669, \notag\\[1.8mm]
    D_1=R_2/R_1=0.999912,&&D_2=0.999535,~~D_3=1.00026, \notag\nd
and by (\ref{9}), $k(n)$ can be determined for $n=1,2,3$ as
\be  k(1)=-0.999823, \quad~k(2)=-1.002791, \quad~k(3)=-0.996851. \notag\nd
As a result, it seems pertinent to regard each of these $k(n)$'s retaining such a uniformly small error from the round figure of $\breve{k}=-1$ to lend it as the rounded value for $k$ with sufficient verifications.  Similarly, we may let each of $z_1(n)=1/R_n$ assume the rounded value of $\breve{z}_1=1.5$ with supporting checks.  We therefore obtain the final result as
\be  \breve{k}=-1, \quad~~ \breve{z}_1=1.5\quad~~\breve{M}=z_1c_0=3, \notag\nd
in resolving the Conjecture for this case with $f(z)$ involving just one simple pole at $z_1=1.5$, of magnitude $M=3$.  Finally, we can draw from the rounded coefficients $\breve{c}_n$'s their specific errors as
\be  \breve{c}_n=c_n+\epsilon_n \quad~\longrightarrow\quad~ \epsilon_1=1.05\times 10^{-4},~~~\epsilon_2=-8.2\times 10^{-5}, ~~~\epsilon_3=2.35\times 10^{-4}, \cdots, \notag\nd
with similar error estimates for $|\breve{D}_n-D_n|,~\ni \breve{D}_n-D_n \rightarrow 0~\forall~\epsilon_n~\rightarrow 0$.  Here we point out that the final solution can be used to re-generate the Fourier coefficients for an ultimate verification of the errors.
\vskip 0.8mm
We also see that this case can suffice to exemplify all the others in this group.  In general, it needs merely several leading terms of series (\ref{5}) to complete the relevant resolution as so illustrated.
\vskip 1mm
\noindent{\bf 4. Advanced Conjecture for the generalized function $f(z)$.}
\vskip 0.8mm
We can extend the original Conjecture to another level of higher interest as follows.
\vskip 0.8mm
\noindent{\bf Advanced Conjecture:} ~{\it The new inverse problem is to have only one of the conjugate functions of $f(z)=u(x,y)+v(x,y)$, say $v$ numerically prescribed on contour $C$ of the unit circle (on which $z=e^{i\theta}, v=\tilde{v}(\theta) ~,-\pi\leq\theta\leq\pi$), other premises being equal, then the new Conjecture asserts that the exact singularity distribution of $f(z)\forall z\in {\cal D}^-$ outside $C$ exists and can be determined.}
\vskip 0.8mm
For resolution, we first separate the series variables into real and imaginary parts as
\begin{subequations}\label{11}
\be c_n=a_n+ib_n, ~~c_ne^{in\theta}&=&(a_n\cos n\theta-b_n\sin n\theta)+i(b_n\cos n\theta+a_n\sin n\theta)\equiv A_n+iB_n, \notag\\
 \tilde{f}(\theta)&=&\sum_{n=0}^N ~c_n e^{in\theta}=\sum_{n=0}~A_n+iB_n=\tilde{u}(\theta)+i\tilde{v}(\theta). \label{11a}\nd
With $\tilde{v}(\theta)$ given, $\tilde{u}(\theta)$ can be derived by applying the circular Hilbert transform (cf Wu[1] Eq(33)):
\be u(\theta)= \hat{H}[v(\phi)]=\frac{\cal P}{2\pi}\int_{-\pi}^\pi v(\phi)\cot \frac{\phi -\theta}{2}d\phi, ~~
     v(\phi) = \hat{H}^{-1}[u(\theta)]=-\frac{\cal P}{2\pi}\int_{-\pi}^\pi u(\theta)\cot \frac{\theta-\phi}{2}d\theta, ~~\label{11b} \nd
which holds for normalized $u$ and $v$ such that $\int_{-\pi}^\pi u(\phi)d\phi =\int_{-\pi}^\pi v(\phi)d\phi=0$.  Hence
\be u(\theta)= \hat{H}[v(\phi)]= \hat{H}[\sin n\phi]= \cos n\theta, \quad v(\phi) = \hat{H}^{-1}[\cos n\theta]= -\sin n\phi \quad~~(n\geq 1). \label{11c} \nd\end{subequations}
Using these relations with the leading term for $n=0$ already known in (\ref{11a}) yields
\be  \tilde{u}(\theta)= \hat{H}[\tilde{v}(\phi)] \quad~\longrightarrow\quad \tilde{f}(\theta)=\tilde{u}(\theta)+i\tilde{v}(\theta) \quad~(-\pi\leq\theta\leq\pi), \label{12} \nd
thus proving that given $\tilde{v}(\theta)$ renders $\tilde{f}(\theta)=\tilde{u}(\theta)+i\tilde{v}(\theta)$ completely prescribed over contour $C$, and hence validating the final solution in (\ref{9}) for the primary resolution to the Advanced Conjecture.
\newpage
\noindent{\bf 5. The Complement Conjecture for the complement function $F(z)$.}
\vskip 0.8mm
The analytical scheme devised above for resolving the Advanced Conjecture can be further extended to exhibit the symmetry in having the roles of domain ${\cal D}^+$ and ${\cal D}^-$ interchanged by converting the Conjecture for the function $f(z)$ reflected into the contour $C$ of unit circle to the Complement Conjecture for the complement function $F(z)$, defined as being regular of order $C^n \forall z\in {\cal D}^-$ and in a neighborhood striding across contour $C$ (cf. Wu[1]). Being so defined, it starts with its integral theorem that
\be  J^-[F(z)]=\oint_{C^-}F(z)dz = -\oint_CF(z)dz=0,  \label{13}\nd
where contour $C^-$ is on the unit circle $|z|=1$, taken in clockwise direction in its own positive sense.  By analogy, the complement function must have a singularity distribution inside contour $C$ (unless $F(z)=$const. on $C$), for which fact we have the analogous conjecture as follows.
\vskip 1mm
\noindent{\bf The Complement Conjecture:}~{\it If the complement function $F(z)$ has on contour $C~(\forall t=e^{i\theta})$ one of its conjugate functions of $F(e^{i\theta})=\hat{U}(\theta)+i\hat{V}(\theta)$, say $\hat{V}$ numerically prescribed, which is equal in value to $\hat{v}(\theta)$ of $f(e^{i\theta})=\hat{u}(\theta)+i\hat{v}(\theta)$ given on $C$, then the new Conjecture asserts that the exact singularity distribution of $F(z)\forall z\in {\cal D}^+$ inside $C$ exists and can be determined.}
\vskip 0.8mm
For resolution, we note that invoking $\hat{V}(\theta)=\hat{v}(\theta)$ dictates that on contour $C~(\forall t=e^{i\theta})$,
\begin{subequations}\label{14}
\be && F(e^{i\theta})=\hat{F}(\theta)=\hat{U}(\theta)+i\hat{V}(\theta)=-\overline{f(e^{i\theta})}=-\hat{u}(\theta)+i\hat{v}(\theta), \label{14a}\\[2mm]
  &&\longrightarrow~~~\hat{U}(\theta)=-\hat{u}(\theta), \quad~ \hat{V}(\theta)=\hat{v}(\theta) \qquad~(-\pi\leq\theta\leq \pi).  \label{14b}\nd\end{subequations}
where $\overline{f(t)}$ denotes the complex conjugate of $f(t)$.  This therefore implies that on $C~(\forall t=e^{i\theta})$,
\begin{subequations}\label{15}
\be  -F(t)=\overline{f(t)}=\overline{f}(e^{-i\theta})= \overline{L}\Sigma_{n=0}(-\bar{z}_1)^{-n}\gamma_n(k) e^{-in\theta}~(=\Sigma_{n=0}\bar{c}_n e^{-in\theta})\qquad (-\pi\leq\theta\leq\pi),\label{15a}\nd
in which the first step follows from (\ref{7}) whereas the second from (\ref{5}), with $L=Mz_1^k$.  This series is convergent and can be continued analytically to all $z=re^{i\theta}$ over the entire $z$-plane, thus giving
\be  -F(z)=\Sigma \bar{c}_n(re^{i\theta})^{-n}=\Sigma \bar{c}_n z^{-n}=\overline{L}\Sigma_{n=0}(-\bar{z}_1)^{-n}\gamma_n(k) z^{-n}=\overline{M}~(\bar{z}_1-z^{-1})^{k},\label{15b}\nd \end{subequations}
where $F(z)$ has its imaginary part equal to that of $f(z)$ on $C$.  This $F(z)$ has a singularity at $z=1/\bar{z}_1$ plus a coupled singularity at $z=0$ (in reflection to that at $z=\infty$ for $f(z)$) which can be made single-valued when $k$ is not an integer with a branch cut between $z=0$ and $z=1/\bar{z}_1$ lying within $C$.
\vskip 0.8mm
Concluding on these extended cases, we note that in addition to the feature that on contour $|z|=1$, $f(z)$ and $F(z)$ have their imaginary parts equal and their real parts only opposite in sign, they further satisfy Cauchy's two integral theorems as expected, i.e.
\begin{subequations}\label{16}
\be   J[f(z)]\equiv\frac{1}{2\pi i}\oint_C \frac{f(t)}{t-z}dt &=& f(z) \qquad~~~ (z\in \overline{{\cal D}^+}={\cal D}^+ +C), \label{16a}\\
    &=& 0   \qquad\qquad~ (z\in \overline{{\cal D}^-}={\cal D}^- +C), \label{16b} \\[1.6mm]
 J^-[F(z)]\equiv\frac{-1}{2\pi i}\oint_C \frac{F(t)}{t-z}dt &=& 0 \qquad\qquad~ (z\in \overline{{\cal D}^+}={\cal D}^+ +C), \label{16c}\\
    &=& F(z)  \qquad~~ (z\in \overline{{\cal D}^-}={\cal D}^- +C). \label{16d} \nd\end{subequations}
These integral formulas are presented specifically to hold for each closed domain in virtue of their corresponding theorems presented in Wu[1].

\vskip 1mm
\noindent{\bf 6. Resolution for function with complex conjugate pairs of singularities.}
\vskip 0.8mm
A case of special interest is for function $f(z)$ having a complex conjugate pair of singularities.  In this case, the series for $f(z)$, in virtue of Schwarz's symmetry, has its coefficients all real.  Here it is apt to first consider a $f(z)$ having a conjugate pair of logarithmic singularities at $z_1=e^{i\alpha}$ and $\overline{z}_1$,
\be f(z)=\log [(1-e^{-i\alpha}z)(1-e^{i\alpha}z)]^{-1/2}=\Sigma_1^\infty~a_nz^n \equiv s(z),\qquad a_n=\frac{1}{n}\cos n\alpha,
\label{17}\nd
$s(z)$ being the {\it series function} defined by the series ($\forall |z|<1$) versus the {\it sum function} $f(z)$ for all $z$.  In case when the series function $s(z)$ is the only data available for locating the singularities of $f(z)$ while $f(z)$ is still undetermined, this particular series can actually provide a useful basic reference.
\vskip 0.6mm
\noindent{\bf Sign pattern criteria of series,} ~by which certain definite relations can be attained between $\alpha$ (hence the argument of $z_1$) and the resulting sign pattern displayed in the series.  For this particular series (\ref{17}), the sign of $a_n$ will depend on the sign of $\cos (n\alpha)$ and will result in an interlacing sign pattern with first $N_1$ terms of positive $a_n$'s, then $N_2$ terms negative, and again $N_3$ terms positive, $\cdots$, so that
\be  -\frac{\pi}{2}<N_1\alpha <\frac{\pi}{2}, ~~~ \frac{\pi}{2}<(N_1+N_2)\alpha <\frac{3\pi}{2}, ~~\cdots,~~ (j-\frac{3}{2})\pi<(N_1+\cdots +N_j)\alpha <(j-\frac{1}{2})\pi, \notag\nd
for $j=1, 2, \cdots$.  In the limit, we actually have just shown a theorem by J. Li[l4](1982) in extending the pioneering work on this subject by M.D. Van Dyke[15](1974).
\vskip 0.5mm
\noindent{\bf Theorem 2. Li's Theorem on Location of a singularity pair of a series.} ~{\it If the power series $\Sigma a_nz^n$, with $a_n\in R, ~\lim |a_n|^{1/n}=1$, has in turn $N_1$ terms positive, $N_2$ terms negative, $N_3$ terms positive, $\cdots$, then the series has a conjugate pair of singularities at $z_1=e^{i\alpha}$ and $\overline{z}_1=e^{-i\alpha}$, where}
\be  \alpha= \lim_{j\rightarrow \infty} j\pi/(\Sigma_{k=1}^j ~N_k). \label{18}\nd
In particular, (\ref{17}) reduces for $\alpha =0$ to $f(z)=-\log (1-z)=\Sigma_1 z^n$, for which $N_1=\infty$, hence (\ref{18}) also yields $\alpha =0$.  Likewise, for $\alpha =\pi$, (\ref{17}) becomes $f(z)=-\log (1+z)=\Sigma_1 (-z)^n$, for which $N_1=N_2=\cdots =1$, then $j\pi/(\Sigma_{k=1}^j ~N_k)=j\pi/j=\pi$ for $j\geq 2$, or as $j\rightarrow \infty$, thus showing the singularity located at $z=-1$ just as that of $f(z)=-\log (1+z)$.  From these two better known cases, we can also infer that $\alpha=\pi/2$, or $\alpha=\pi/3$, or $\alpha=\pi/N$ if the signs of $a_n$'s interlace in pairs, or in triplets, or in $N$-tuples, respectively.
\vskip 0.6mm
In general, the angular positions of all logarithmic singularity pairs can be determined as a sum of the individual pairs.
Further generalizing Li's Theorem, we can readily resolve the singularities of products of multiple pairs of algebraic conjugate functions such as $f(z)= (z-z_1)^{-k}(z-\overline{z}_1)^{-k} ~(z_1=e^{i\alpha})$, since we can then set
\be g(z)=\log f(z) =-k \log \{(1-z_1 z)(1-\overline{z}_1 z)\}= 2k ~\Sigma_1^\infty ~n^{-1} z^n \cos n\alpha, \notag\nd
which is again in the form of (\ref{17}), and therefore giving the singularities of $f(z)$ as that of $f(z)=\exp \{g(z)\}$.  The same logical argument holds for other types of such functions.
\vskip 1mm
\noindent{\bf 7. Resolution to the Conjecture for $f(z)$ having multiple singularities.}
\vskip 0.8mm
For more general cases, we consider next the case of $N_s\geq 2$ for singularities of two kinds, one being those of equal order and arbitrary locations whereas the other for those of arbitrary orders and locations, to begin with the case of $N_s=2$ for possible extension to multiple singularities.
\vskip 1mm
\newpage
\noindent{\bf 7.1. ~Resolution for $f(z)$ having two singularities of equal order.}
\vskip 0.8mm
It is of interest to consider $f(z)$ having multiple singularities of equal order, such as that found in the literature on water waves.  In the case of $N_s=2$ for two singularities of equal order, we have
\be  f(z)&=&M_1(z_1-z)^k +M_2(z_2-z)^k =\Sigma_{n=0}^\infty(L_1 z_1^{-n}+L_2 ~z_2^{-n})\gamma_n(k)(-z)^n, \label{19}\nd
where $L_i=M_i z_i^k ~(i=1,2)$ and $\gamma_n(k)$ is given in (\ref{7a}).  Similar to (\ref{5})-(\ref{9}), matching series (\ref{19}) for $z=e^{i\theta}$ termwise with the prescribed series in (\ref{5}) yields
\be (L_1 ~z_1^{-n}+L_2~z_2^{-n})\gamma_n=(-1)^n c_n \qquad~(n=0,1,2,\cdots), \label{20}\nd
These are the set of transcendental algebraic equations for solving the five base parameters in total, $(M_1,M_2,z_1,z_2, k)$ and for verifying the consistency for this case.
\vskip 0.5mm
First, we eliminate the terms with $L_2$ by operation for $\{(k-n)/(n+1)z_2\}c_n+c_{n+1}$, giving
\begin{subequations}\label{21}
\be  L_1(z_2^{-1}-z_1^{-1})\gamma_{n+1}z_1^{-n}=(-1)^n \tilde{c}_{n+1}, \quad~~\tilde{c}_{n+1}=\frac{k-n}{(n+1)z_2}c_n+c_{n+1}. \label{21a}\nd
Next, $L_1$ can now be eliminated by taking the ratio $\tilde{c}_{n+1}/\tilde{c}_n ~(n=1, 2,\cdots)$, yielding
\be \frac{\tilde{c}_{n+1}}{\tilde{c}_n}=-\frac{\gamma_{n+1}}{\gamma_{n}z_1}=\frac{n-k}{n+1}\frac{1}{z_1}=\frac{n}{n+1}H_n(k,z_2) \quad~~ H_n(k,z)=\frac{(n+1)R_nz-(n-k)}{nz-(n-1-k)/R_{n-1}}, \label{21b}\nd
which can clearly be written in symmetry between the two singularities in equality as
\be  G_n(k,z_1)=G_n(k,z_2), \quad~~ G_n(k,z)=H_n(k,z)/z \qquad~(n=1,2,\cdots). \label{21c}\nd\end{subequations}
Finally, expanding out (\ref{21b}) yields the basic equation for this case as
\be  (n+1)R_nz_1z_2-(n-k)\{(z_1+z_2) -(n-k-1)/(nR_{n-1})\}=0 \quad~~ (n=1,2,\cdots), \label{22}\nd
which also exhibits the symmetry between the two singularities interchangeable in designation.
\vskip 0.6mm
To resolve the three unknown parameters $(k,z_1,z_2)$ in (\ref{22}), we may take four leading equations of (\ref{22}) to eliminate $(z_1+z_2)$ and $z_1z_2$ in two steps, yielding one equation for $k$ ($n=1,2,\cdots$) as
\begin{subequations}\label{23}
\small
\be \frac{K(-1,k)}{R_{n-1}}\left(\frac{R_{n+1}}{K(1,k)}-\frac{R_{n+2}}{K(2,k)}\right)+ \frac{K(1,k)}{R_{n+1}}\frac{R_n}{K(0,k)} +\frac{K(0,k)}{R_n}\frac{R_{n+2}}{K(2,k)}=2, \quad~K(j,k)=\frac{n+j-k}{n+1+j},\normalsize \label{23a}\nd
which is a cubic equation for $k$ (after dismissing the complex conjugate roots by confirmation).  With $k$ so determined, ($z_1z_2$) can be deduced by eliminating $z_1+z_2$ from (\ref{22}), giving for $n=1,2,\cdots$,
\be  z_1z_2=\left(\frac{n-k}{(n+1)R_n}-\frac{n-k-1}{R_{n-1}}\right)/\left(\frac{n+1}{n-k}R_n-\frac{n+2}{(n+1)-k}R_{n+1}\right)\equiv A(n,k)  \label{23b}\nd
for $n=1$, say, with which $z_1+z_2 ~(=2B(n,k)$, say) then follows from (\ref{22}) for $n=1$.  Finally, $z_1$ and $z_2$ are obtained by combining $z_1z_2=A(n,k)$ and $z_1+z_2=2B(n,k)$ to give the quadratic equation,
\be  z_1^2-2Bz_1+A=0 ~~\longrightarrow ~~~z_1=B+\sqrt{B^2-A},~~~z_2=B-\sqrt{B^2-A}. \label{23c}\nd\end{subequations}
At last, the solution is completed with the singularity magnitudes found from (\ref{21a}) for $n=0$ and verified for its feature of consistency using (\ref{19}) for $n\geq 5$.  Or in another approach, we may use ($M_1,k,z_1$) Just obtained to provide its contribution, say ${c}'_n ~(n=0,1,\cdots)$ to differ from the original Fourier coefficients $c_n$ by $c_n-{c}'_n$ as the new data for the other singularity ($M_2,k,z_2$) by direct application of (\ref{8}) for verifying the solution.  Our resolution to the Conjecture is then accomplished for this case of $f(z)$ having two singularities of equal order.  When there are more than two singularities of equal order, this method can again provide a scheme for solution by induction.
\vskip 1mm
\noindent{\bf 7.2. ~Resolution for $f(z)$ having two individual singularities.}
\vskip 0.8mm
For $f(z)$ having two singularities of arbitrary orders and locations, we have
\be  f(z)&=&M_1(z_1-z)^{k_1} +M_2(z_2-z)^{k_2} =\Sigma_{n=0}^\infty(L_1 z_1^{-n}\gamma_n(k_1)+L_2 ~z_2^{-n}\gamma_n(k_2))(-z)^n, \label{24}\nd
where $L_i=M_i z_i^{k_i} ~(i=1,2)$.  Again, matching series (\ref{24}) with series in (\ref{5}) yields for this case
\begin{subequations}\label{25}
\be L_1 ~z_1^{-n}\gamma_n(k_1)+L_2~z_2^{-n}\gamma_n(k_2)=(-1)^n c_n \qquad~(n=0,1,2,\cdots), \label{25a}\nd
These are the set of equations for solving $(M_i,z_i, k_i,~i=1,2)$ and for verifications.  Elimination of $L_1$ and $L_2$ can be carried out in close analogy with that for (\ref{19}), yielding the basic equation as
\be \frac{z_1-z_2}{(k_2-n+1)z_1-(k_1-n+1)z_2}= 1+\frac{nz_1}{k_1-n+1}\cdot\frac{(k_2-n) +(n+1)R_nz_2}{n z_2-(n-k_2-1)/R_{n-1}},  \label{25b}\nd\end{subequations}
which reduces to (\ref{22}) for $k_1=k_2=k$.  Pursuing the image symmetry that (\ref{25b}) is invariant if $(z_1, k_1)$ and $(z_2, k_2)$ are interchanged between the two independent singularities, we obtain
\be  G_n(k_1, z_1) =G_n(k_2, z_2), \quad~~ G_n(k,z)=H_n(k,z)/z \qquad~(n=1,2,\cdots), \label{26}\nd
where $H(k,z)$ is given in (\ref{21b}).
\vskip 0.8mm
Now, the four remaining unknowns, $k_1,k_2,z_1,z_2$ can be solved by taking four equations of (\ref{25b}) or (\ref{26}), for $n=1,2,3,4$, with the rest serving for verifying the solution.
\vskip 0.8mm
\noindent{\bf Example 2.} ~Here we have the coefficients in numerics for series (\ref{5}) as
\be c_1=4, \quad c_2=1.5, ~~~c_3=0.583333, ~~~c_4=0.236111, ~~~c_5=0.099537, ~~~c_6=0.043596, \cdots,~\longrightarrow ~~&&\notag\\[1mm]
 R_1=0.375, ~~~~R_2=0.388889, ~~~~R_3=0.404762, ~~~~R_4=0.421569, ~~~~R_5=0.437985, \cdots,\qquad~~~&&\notag\nd
Not knowing the number of singularities, we first test for $N_s=1$.  Then, like in Example 1, we find
\be  k(1)=-1.074074,~~~k(2)=-0.781818,~~~k(3)=-0.572727, ~\cdots, \notag\nd
for the singularity key index $k(n)$ of the few leading orders.  Their lack of any feature of consistency evidently indicates the singularity distribution being not in the single singularity group.
\vskip 0.8mm
So next, we test out the group of two separated singularities for this $f(z)$.  By substituting the above $R_n$'s in the four leading equations of (\ref{25b}) to solve for $(k_1,k_2,z_1,z_2)$, adopting successive elimination or any other adequate algorithms, we obtain, after some algebra, the final solution as
\be  k_1=k_2=-1, ~~~z_1=2,~~~z_2=3,~~~M_1=-2,~~~M_2=-9, \notag\nd
in which $M_1, M_2$ follow from (\ref{25a}) for $n=0$ and $n=1$, with the consistency exhibited as
\be  G_1(z_1,k)=G_1(z_2,k)=-\frac{1}{3},~~~G_2(z_1,k)=G_2(z_2,k)=-\frac{1}{2},~~~G_3(z_1,k)=G_3(z_2,k)=-\frac{2}{3}, \notag\nd
followed by $G_4(z_1,k)=G_4(z_2,k)=-5/6$, each with a relative error of $10^{-8}$, and likewise for $n\geq 5$.
\vskip 0.8mm
\noindent{\bf Example 3.} ~As the two singularities in Example 2 turn out to be of equal order, we can use the same coefficients $c_n$ for a comparative study.  Substituting the $R_n$'s of Example 2 in (\ref{23a}) yields $k=-1$, with the other two complex conjugate roots dismissed by verification.  With $k=-1$, (\ref{23b}) for $n=1$ produces $z_1z_2=6$ and then (\ref{22}) gives $z_1+z_2=5$.  Therefore by (\ref{23c}), we obtain the solution
\be  k=-1, \quad~ z_1=2, \quad~z_2=3, \notag\nd
in agreement with Example 2, here also with relative errors of O($10^{-8})$.  In comparison, the algebra involved here is noticeably simpler than the former case, especially for conducting verifications.
\vskip 0.8mm
When more than two singularities are encountered, this case for two singularities can be expected to provide a mathematical scheme for resolution by induction, i.e. by eliminating the magnitude of one singularity at a time till all the magnitude variables are eliminated for solving the locations of the singularities and their key order indices by complex algebra.  However, the computation along this approach is still not simple, whereas a more general approach can now be addressed next.
\vskip 1mm
\noindent{\bf 8. Resolution to the Conjecture involving multiple singularities by induction.} ~
\vskip 0.8mm
For this general case, we reconsider the group of functions having multiple singularities, all lying outside the open domain ${\cal D}^+$ bounded by the unit circle $C$, and situated all at different radial distances from the origin (as the case of singularities situated on a concentric circle beyond $C$ has been discussed in \S 6).  To begin with, we first consider $f(z)$ having two singularities at $z_1$ and $z_2$, both real and positive ~$(1<z_1<z_2)$, a special case which is apt for illustrating the underlying basic principle more easily, and can be extended to functions having more singularities situated more freely.
\vskip 0.8mm
To have the singularities of $f(z)$ determined by induction, we first consider the direct problem.
\vskip 0.8mm
\noindent{\bf 8.1. The direct problem.} ~{\it The direct problem has the singularities all given in closed analytical form as a basis for developing a general method to determine the singularities by induction so that the method can readily be adapted to resolving the inverse problem as stated by the Conjecture.}
\vskip 0.6mm
Thus, we take, parallel to (\ref{24}), here with $1<z_1<z_2$,
\be  f(z)=f_1(z)+f_2(z)= \Sigma_{j=1}^2M_j(1-z/z_j)^{k_j}=\Sigma_{n=0}^\infty\Sigma_j~M_j\gamma_n(k_j)(-z/z_j)^n=\Sigma_n \hat{c}_n z^n, ~~~~&&\notag\\[1.2mm]
 \hat{c}_n =M_1\gamma_n(k_1)(-z_1)^{-n}\left\{1+m\frac{\gamma_n(k_2)}{\gamma_n(k_1)}\lambda^n\right\} \quad(m=M_2/M_1, ~0<\lambda=z_1/z_2<1), &&\label{27}\nd
where $\hat{c}_n$ derives its value from the analytic formulas given for this direct problem to give
\be R_n=\frac{\hat{c}_{n+1}}{\hat{c}_n}=\frac{1-k_1/n}{1+1/n}\cdot\frac{Q_n}{z_1}, \quad~~Q_n=\frac{1+m(\gamma_{n+1}(k_2)/\gamma_{n+1}(k_1))\lambda^{n+1}}{1+m(\gamma_n(k_2)/\gamma_n(k_1))\lambda^n}. \label{28}\nd
For $0<\lambda<1$, $\lambda^n\rightarrow 0$ and $Q_n\rightarrow 1$ as $n\rightarrow\infty$.  Hence, for given $\lambda<1$, there exists $N \ni\forall n\geq N,~|1-Q_n|<\epsilon$.  Thus, as $\epsilon\rightarrow 0$, $R_n$ tends to that for a single singularity, which can then be determined in this asymptotic limit when the effects due to the singularity farther away fall off.
\vskip 0.8mm
\noindent{\bf Domb-Sykes plot.} ~Proceeding for this asymptotic solution, it is apt to make the Domb-Sykes[11] plot of $R_n$ versus $1/n$ as a function $R_n(1/n)$, just as so conspicuously exhibited by the functional relation in (\ref{28}), with $(k_1,z_1,Q_n(\lambda))$ as the primary parameters.  The basic principle and the main features of the solution can be illustrated below.
\vskip 0.8mm
\noindent{\bf Example 4. Function $f(z)$ having two logarithmic singularities.} ~Of basic interest is that
\be f(z)=\log\{(1-z/z_1)(1-z/z_2)^2\}=\log (1-z/z_1)+2\log (1-z/z_2) \quad~~(z_1=1.1, ~z_1<z_2), \label{29}\nd
Here, we have $f(z)=f_1(z)+f_2(z), ~k_1=k_2=0, ~m=2, ~z_1=1.1, ~\lambda=z_1/z_2 ~(0<\lambda<1)$, hence
\be  R_n=\frac{n}{n+1}\cdot\frac{Q_n}{z_1}, \qquad~ Q_n=\frac{1+2\lambda^{n+1}}{1+2\lambda^n}\qquad~(n=1,2,\cdots). \label{30}\nd
This function $R_n(1/n)$ has been computed by employing Mathematica 7 for three values of $\lambda=0.125, ~0.4, ~0.8,$ over the range $n=1, 2,\cdots, 100$, with the resulting $R(1/n)$ plotted in three lines of sequence of points for each $n$ as shown in Fig. 1, with the $\Diamond$ indicating the points for $\lambda=0.125$, the $\Box$ for $\lambda=0.4$, and the $\bullet$ for $\lambda=0.8$.  Here, for $n\geq N \simeq 20$, the three lines of dotted points apparently merge into one straight line, ending at $1/n=0.01$ (or $n=100$).  This line is further extrapolated by being fitted to $R=(p_o+p_1/n)/(1+q_1/n)$ to reach the $R$-axis intercept at $R(0)=0.9091$, of slope $R'(0)=-1$ which gives $k=0$, thus implying $f_1(z)$ logarithmic and further $z_1$ being located at $z_1=1/R(0)=1.1$, with almost no error.  These are the two important data found from the plot.
\begin{figure}[!h]
\begin{center}
\includegraphics[width=0.6\textwidth]{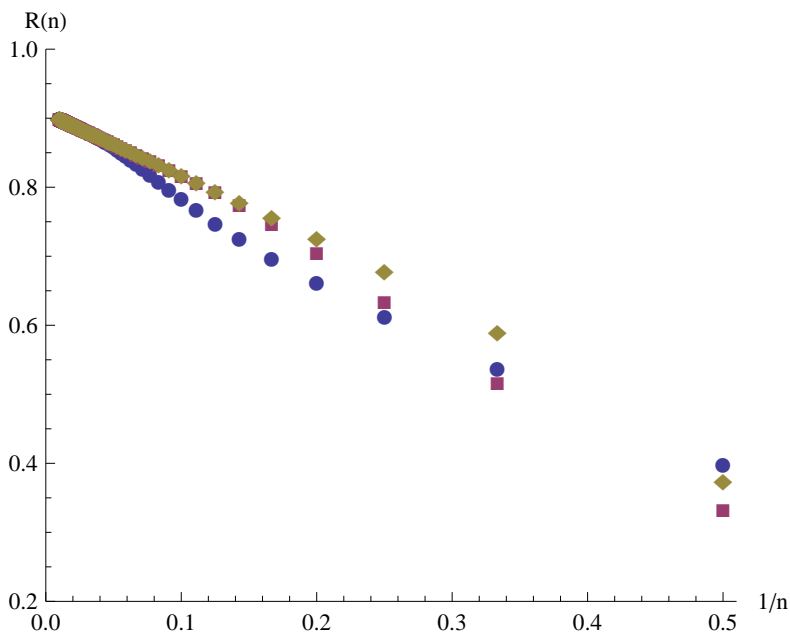}
\caption{\small Domb-Sykes plot of $R_n(1/n)$ for three values of
singularity radial-distance ratio $\lambda=z_1/z_2$, shown by points
$\Diamond$ for $\lambda=0.125$, by $\Box$ for $\lambda=0.4$, by
$\bullet$ for $\lambda=0.8; n=1,\cdots, 100.$}
\end{center}
\end{figure}
It is of interest to note that for $n\leq N\simeq 20$, the three dotted curves have the $\Diamond$ line (for $z_2=8 z_1$) staying the closest to the asymptotic line, $L_a: R_a=R(0)+R'(0)/n$, whereas the $\Box$ line (for $z_2=2.5 z_1$) and the $\bullet$ line (for $z_2=1.25 z_1$) bifurcate increasingly more from the straight $L_a$ line with decreasing radial distances.  This signifies that the interaction between two point singularities increases with decreasing distances apart.  This is another key feature demonstrated by the Domb-Sykes plot.
\vskip 0.8mm
At this point, we note that the number of singularities, $N_s=2$, and $k=0$ being both known in this direct problem, we can directly solve for the other singularity by taking the first two $R_1, ~R_2$ to deduce from their formulas to attain  $m=2$ and the values of $z_2$ for the three cases, and their magnitudes from $c_o$ to complete the solution to this direct problem.
\vskip 0.8mm
\noindent{\bf 8.2. The inverse problem.} ~Proceeding to resolve the Conjecture, we now take up the inverse problem which has only the Fourier coefficients $c_n$ of (\ref{5}) prescribed in numerics, ordinarily by the solution to a problem solved by a numerical scheme, but here assumed to take the value $\hat{c}_n$ from the above direct problem for a single case for $\lambda=0.4$, say.  With $c_n=\hat{c}_n$ known, we plot the coefficient ratio $R_n=c_{n+1}/c_n$ to attain in this case again the same result as that shown in Fig. 1, thereby having $(k_1=0, ~z_1=1.1)$ determined for the first singularity of $f_1(z)$ (lying nearest to the unit circle $C$) as
\be  f_1(z)= M_1\log (1-z/z_1)= c'_o -\Sigma_{n=1}~c'_n z^n \quad~~(c'_o=0, ~c'_n=M_1z_1^{-n}/n, ~n\geq 1), \label{31}\nd
where $c'_o=c_o=0, ~M_1=1$ by deduction from $c'_n=c_n$ asymptotically, and hence $f_1(z)$ is all known.
\vskip 0.8mm
To this end, it is vital to note that the dotted line (here with $\Box$) being curved indicates that $f(z)$ has more singularities, possibly not just one more, to complete for resolving the inverse problem.
\vskip 0.8mm
\noindent{\bf 8.3. Resolving the inverse problem by induction.}~To continue, we rid $f(z)$ of $f_1(z)$, giving
\be f_2(z)=f(z)-f_1(z)=\Sigma_{0}~c''_nz^n,\quad~c''_n=c_n-c'_n ~~~\longrightarrow~~~ R''_n=c''_{n+1}/c''_n \quad~~ (n=1,2,\cdots). \label{32}\nd
With $c''_n$ here all known, we can carry out the Domb-Sykes plot for $R''_n(1/n)$ to determine the primary parameters $(k_2,z_2)$ for $f_2(z)$.  In this particular case for $f(z)$ of (\ref{29}), we will find from the $R''_n(1/n)$ plot that the entire dotted line is straight, of slope $dR''/d(1/n)=-1$, and the $R''$-intercept at $1/z_2=0.36364$ for $z_2= 2.75=2.5 z_1,$ as the last singularity of $f(z)$, therefore completing this inverse problem.
\vskip 0.8mm
In general, for generic cases, the $R''(1/n)$ line may appear again curved, in which case we repeat the analogous plot for $R^{(3)}(1/n))$ for $f_3(z), \cdots$, until $R^{(\ell)}(1/n))$ for $f_\ell(z)$ bearing out a straight line in the plot for the very last singularity of $f_\ell(z)$ to complete the inverse problem.  This then constitutes the general method of resolving the Conjecture by induction.
\vskip 1mm
\noindent{\bf 9. Resolution to the Conjecture involving essential singularities.} ~
\vskip 0.8mm
To enhance addressing this issue, we apply the generalized Hilbert transforms (cf. Wu[1]).
\vskip 0.8mm
\noindent{\bf 9.1. Classical Hilbert transform.} ~First, we consider a class of analytic $C^1$ function $f(z)=u(x,y) + iv(x,y)$ which is regular in the upper half $z$-plane for $Im~z\geq 0$, and vanishes as $|z|\rightarrow \infty$ uniformly in $0\leq \arg z\leq \pi$, with the upper-half (or lower-half) $z$-plane designated as the open domain ${\cal D}^+$ (or ${\cal D}^-$), bounded by the $x$-axis ($y=0$) and the point $z=\infty$.  Then we have the Hilbert transform between $u(x)$ and $v(x)$ along the $x$-axis given by
\be  u(\xi)=H[v(x)]=\frac{1}{\pi}{\cal P}\int_{-\infty}^\infty\frac{v(x)dx}{x-\xi}, \qquad
 v(x)=H^{-1}[u(\xi)]= \frac{-1}{\pi}{\cal P}\int_{-\infty}^\infty\frac{u(\xi)d\xi}{\xi-x},  \label{33} \nd
with $H$ denoting the Hilbert transform and $H^{-1}$ the inverse transform.  We note that as a point $z$ tends to $(x,0)$ from ${\cal D}^+$ (or ${\cal D}^-$), $f(z)$ in (\ref{1a})-(\ref{1b}) tends in the limit to $f^+(x)=f(x),~f^-(x)=0, ~f(x)=u(x)+iv(x), ~u(x)$ and $v(x)$ being given by (\ref{33}), all as a special case of (\ref{3}) for this particular geometry.
\vskip 0.8mm
\noindent{\bf Example 5. Hilbert transform of sinusoidal functions.} ~For $v(x)=\sin x$, we have
\be  u(\xi)=H[\sin(x)]=\frac{1}{\pi}{\cal P}\int_{-\infty}^\infty\frac{\sin x~dx}{x-\xi}= \cos (\xi). \label{34} \nd
Hence, $f(x)=u(x)+iv(x)= e^{ix}$ on the $x$-axis, which can be continued into the entire $z$-plane, giving
\be  f(x)=u(x)+iv(x)= e^{ix}\quad~~\longrightarrow\quad~~ f(z)=e^{iz}\quad~~(\forall ~z \in 0\leq |z|<\infty).\label{34} \nd
Clearly, this function, $e^{iz}=e^{ix-y}$ is regular in the upper-half $z$-plane, vanishing exponentially as $y\rightarrow +\infty$, but becomes essentially singular as $y\rightarrow -\infty$ in domain ${\cal D}^-$ (the lower-half $z$-plane).
\vskip 0.8mm
\noindent{\bf 9.2. Complement Hilbert transform.} ~On the contrary, if complement function $F(z)=U(x,y)+iV(x,y)$ is $C^1$ regular in the open domain ${\cal D}^- \{z|Im ~z \leq 0\}$, then applying the general formula (\ref{3}) produces for this particular geometry in ${\cal D}^- $ the complement Hilbert transform between $U(x)$ and $V(x)$ along the $x$-axis the equations (cf. Wu[1]) as
\be  U(\xi)=\overline{H}[V(x)]=\frac{-1}{\pi}{\cal P}\int_{-\infty}^\infty\frac{V(x)dx}{x-\xi}, \quad~~
 V(x)=\overline{H}^{-1}[U(\xi)]= \frac{1}{\pi}{\cal P}\int_{-\infty}^\infty\frac{U(\xi)d\xi}{\xi-x}, \label{36}\nd
with $\overline{H}$ denoting the complement Hilbert transform and $\overline{H}^{-1}$ the inverse transform.  They can be employed for evaluating transforms of various complement functions.
\vskip 0.8mm
\noindent{\bf Example 6. Complement Hilbert transform of sinusoidal functions.} ~If for the complement function $F(x)=U(x)+iV(x)$ on the $x$-axis, $V(x)=\sin x$ (equal in value to $v(x)=\sin x$ in (34), then by complement Hilbert transform (\ref{36}), we obtain
\be  U(\xi)=\overline{H}[\sin(x)]=-\frac{1}{\pi}{\cal P}\int_{-\infty}^\infty\frac{\sin x~dx}{x-\xi}= -\cos (\xi).\label{37} \nd
Hence by analytical continuation, we have
\be  F(x)=U(x)+iV(x)= -e^{-ix}\quad~~\longrightarrow\quad~~ F(z)=-e^{-iz}\quad~~(\forall ~z \in 0\leq |z|<\infty).\label{38} \nd
Obviously, $F(z)=-e^{-ix+y}$ is regular in ${\cal D}^- $, and has an essential singularity in ${\cal D}^+$ at $z=\infty.$  Therefore, one equal value of one conjugate function, e.g. $v(x)=V(x)=\sin x$, produces an analytic function $f(z)$ regular in ${\cal D}^+$ and another $F(z)$ regular in ${\cal D}^-$, whilst both having an essential singularity in each of their complement domain, respectively, in splendid analogy with the original Conjecture in \S 1 and the Complement Conjecture in \S5, conjointly resulting in (\ref{16}) as their joint roots.

\vskip 1mm
\noindent{\bf 10. Discussion and conclusion.} ~
\vskip 0.8mm
In this study we have illustrated how the solutions can be obtained to the inverse problem for the generalized function $f(z)$ having total number of $N_s=1,2,\cdots$ singularities in resolving the Conjecture.  The simpler cases (i)-(ii) for $f(z)$ having only one or two singularities are seen to have provided the study a sound base by exhibiting a clear structure of all solutions to the inverse problem.  The general method developed here in \S8 for resolution to the Conjecture by induction is suitable for singularities located at differing radial distances so that the one closest to the origin can be singled out by applying the Domb-Sykes plot for an asymptotic solution, leaving the remaining singularities to be similarly resolved one at a time.  This is complemented by the case of singularities lying on a concentric circle by the analytical method expounded in \S7 and by another approach utilizing the series sign pattern criteria that has a longer history with the pioneering works of Van Dyke[15] and Li[14].  The success of these methods and the high accuracy of the general results are seen to stem from the central role played by the power series in complex form composed of orthogonal terms.  These methods can now be employed for application and further development to all pertinent problems in mathematics, applied science, and mathematical physics.
\vskip 0.8mm
As mentioned in foregoing occasions, there are at least two major areas in analysis of holomorphic functions where the existing challenges can be benefited from the resolution of the Conjecture.  The long and rich history of the fully nonlinear theory of dispersive water waves can serve as a splendid representative for all the disciplines in which the key advances may require a final resolution of their principal theory of perturbation expansions.  In various other circumstances, new grasp of the exterior singularities in domain complement to the solution region of the primary variables may furnish crucial data base to ascertain the certitude of solution validity and profound comprehension of the target phenomena.  It may also cast light on the course of expounding the key mechanism underlying instabilities of a nonlinear system at an initial stationary state.
\vskip 0.8mm
In concluding, it is reasonable to say that this is just a beginning.  New issues may arise to require answers.  The exterior singularities of analytic functions may present a variety of new interests, e.g. coalescence of singularities with varying parameters, productions of essential singularities, etc.  We may indeed expect to see more experience accumulating to enrich more of this new field.
\vskip 1.8mm
\noindent{\bf Acknowledgment.} ~I wish to thank very kindly Prof. Joe Keller, Prof. Lu Ting and Prof. John C.K. Chu for their interest in resolving the Conjecture enunciated in a previous paper[1], and especially Prof. Thomas Y. Hou for his very careful scrutiny of the present article with interesting queries.  I would also like to thank Yue Yang for his assistance in preparing the figure, and I am further highly appreciative for the gracious encouragement from Dr. Chinhua S. Wu and the American-Chinese Scholarship Foundation.
\vskip 1.8mm
\noindent{\bf References.}
\vskip 0.8mm
\noindent [1] Wu, Th.Y. On uniform continuity of Cauchy's function and uniform convergence of Cauchy's integral \\
\indent formula with applications.  arXiv:0710.5790v2 [math.CV]29 Dec 2007, pp.1-21 (2007).
\vskip 0.8mm
\noindent [2] Stokes, G.G. On the theory of oscillatory waves. {\it Trans. Cambridge Phil. Soc.} {\bf 8}, 441-455 (1847).\\  \indent In {\it Math. Phys. Papers,} Vol 1, 197-229. Cambridge U. Press, (1880).
\vskip 0.8mm
\noindent [3] Longuet-Higgins, M.S. and Fox, M.J.H. Asymptotic theory for the almost-highest solitary wave. \\
\indent {\it J. Fluid Mech.} {\bf 317}, 1-19 (1996).
\vskip 0.8mm
\noindent [4] Williams, J.M. Limiting gravity waves in water of finite depth. {\it Phil. Trans. R. Soc. Lond.}{\bf A 302}\\ 
\indent 139-188 (1981).
\vskip 0.8mm
\noindent [5] Wu, Th.Y., Kao, J. \& Zhang, J.E. A unified intrinsic functional expansion theory for solitary waves.\\ \indent {\it Acta Mech. Sinica}  {\bf 21}, 1-15 (2005).
\vskip 0.8mm
\noindent [6] Qu, W. Studies on nonlinear dispersive water waves. Ph.D. Thesis, California Institute of Technology.\\
\indent (2000).
\vskip 0.8mm
\noindent [7] Wu, Th.Y., Wang, X. \& Qu, Q. On solitary waves. Part 2. A unified perturbation theory for \\ 
\indent higher-order waves.  {\it Acta Mech. Sinica}  {\bf 21}, 515-530 (2005).
\vskip 0.8mm
\noindent [8] Wang, X. \& Wu, Th.Y. Integral convergence of the higher-order theory for solitary waves. \\
\indent {\it Physics Letters A.} {\bf 350}, 44-50 (2006).
\vskip 0.8mm
\noindent [9] Grant, M.A. The singularity at the crest of a finite amplitude progressive Stokes wave. {\it J. Fluid \\
\indent Mech.} {\bf 59}(2), 257-262 (1973).
\vskip 0.8mm
\noindent [10] Schwartz, L.M. Computer extension and analytic continuation of Stokes expansion of gravity waves. \\
\indent {\it J. Fluid Mech.} {\bf 62}(3), 553-578 (1974).
\vskip 0.8mm
\noindent [11] Domb, C \& Sykes, M.F. {\it Proc. R. Soc. Lond.} {\bf A240} 214-218 (1957).
\vskip 0.8mm
\noindent [12] Tanveer, S. Singularities in water waves and Rayleigh-Taylor instability. {\it Proc. R. Soc. Lond.}\\
\indent {\bf A~435} 137-158 (1991).
\vskip 0.8mm
\noindent [13] Longuet-Higgins, M. S. \& Fox, M.J.H. Theory of the almost-highest wave. Part 2. Matching and \\ 
\indent analytic extension. {\it J. Fluid Mech.} {\bf 85} 769-786 (1978). 
\vskip 0.8mm
\noindent [14] Li, J. {\it Scientia Sinica} ~{\bf A~25}(6) 593-600 (1982).
\vskip 0.8mm
\noindent [15] Van Dyke, M.D. {\it Q. J. Mech. Appl. Math.} ~{\bf 27} 423-440 (1974).

\end{document}